\documentclass[reqno]{amsart}
\usepackage{amsmath}%
\usepackage{amsfonts}%
\usepackage{amssymb}%
\usepackage{graphicx}
\def \P  {{\rm {I\hspace{-.03in}P}}}

\def \pgd {LDP }

\def \cvloi {\stackrel {w}{\rightarrow}}

\def \Cb {C_{b}(\Sigma \times \Sigma)}

\newtheorem{theorem}{Theorem}
\theoremstyle{plain}

\newtheorem{lemma}{Lemma}

\numberwithin{equation}{section}
\begin{document}
\title[LD for symmetrised empirical measures]{Large deviations for symmetrised empirical measures}
\author{Jos\'e Trashorras}
\thanks{Universit\'e Paris-Dauphine, Ceremade, Place du Mar\'echal de Lattre de Tassigny 75775 Paris Cedex 16 France. Email: {\sf xose@ceremade.dauphine.fr}}
\date{Mai 29, 2007}

\begin{abstract}
In this paper we prove a Large Deviation Principle for the sequence of symmetrised empirical measures 
$\frac{1}{n} \sum_{i=1}^{n} \delta_{(X^n_i,X^n_{\sigma_n(i)})}$ where $\sigma_n$ is a 
random permutation and $((X_i^n)_{1 \leq i \leq n})_{n \geq 1}$ is a triangular array of random variables 
with suitable properties. As an application we show how this result allows to improve the Large Deviation 
Principles for symmetrised initial-terminal conditions bridge processes recently established by 
Adams, Dorlas and K\"{o}nig.  
\end{abstract}
\maketitle


\section{Introduction and results}


We say that a sequence of Borel probability measures $(P^{n})_{n \geq 1}$ on a topological space 
$\mathcal{Y}$ obeys a Large Deviation Principle (hereafter abbreviated LDP)
with rate function $I$ if $I$ is a non-negative, lower semi-continuous function defined on $\mathcal{Y}$ 
such that

$$ 
- \inf_{y \in A^{o}} I(y) \leq \liminf_{n \rightarrow \infty} \frac{1}{n} \log P^{n}(A)
     \leq \limsup_{n \rightarrow \infty} \frac{1}{n} \log P^{n}(A) \leq -
     \inf_{y \in \bar{A}} I(y)
$$

\noindent
for any measurable set $A \subset \mathcal{Y}$, whose interior is denoted by $A^{o}$ and closure by $\bar{A}$. 
If the level sets $\{y : I(y) \leq \alpha \}$ are compact for every $\alpha < \infty$, $I$ is called 
a good rate function. With a slight abuse of language we say that a sequence of random variables obeys an 
\pgd when the sequence of measures induced by these random variables obeys an
LDP. For a background on the theory of large  deviations, see Dembo and Zeitouni \cite{dz98} 
and references therein.

\bigskip

\bigskip

Our framework is the following: Let $(\Sigma,d)$ be a Polish 
space and $((X_i^n)_{1 \leq i \leq n})_{n \geq 1}$ be a triangular array of $\Sigma-$valued random variables defined on a 
probability space $(E,\mathcal{E},\mathbb{P})$. Let $(\sigma_n)_{n \geq 1}$ be a sequence of random 
variables defined on the same probability space and such that for every $n \geq 1$ the distribution of $\sigma_n$ 
is uniform over the set $\mathfrak{S}_n$ of all permutations of $\{1,\dots,n\}$. We further assume that 
for every $n \geq 1$ the vector $(X_1^n,\dots,X_n^ n)$ is independent from $\sigma_n$. Our purpose is to 
investigate the LD properties of the sequence of random measures

\begin{equation}
\mathcal{L}^n = \frac{1}{n} \sum_{i=1}^{n} \delta_{(X_i^n,X_{\sigma_n(i)}^n)}
\label{sym}
\end{equation}
which we call symmetrised empirical measures.

Our interest for these measures comes from the recent work of Adams, Bru, Dorlas and K\"{o}nig \cite{a06b,a06a,a06d,a06c}
on the modeling and rigorous analysis of the Bose-Einstein condensation at positive temperature. Indeed, symmetrised 
measures are key ingredients in the modeling of Boson systems i.e. quantum particle systems whose wave function 
are invariant under permutation of the single-particle variables. For such models the thermodynamic limit of the free 
energy is obtained by computing the limit of the properly rescaled trace of the restriction of the Boltzmann factor
$e^{- \beta \mathcal{H}_n}$ to the subspace of symmetric wave functions. In \cite{a06d} and \cite{a06c} this is achieved 
for simple models by combining LDPs for bridge processes with symmetrised initial-terminal conditions
with Feynman-Kac formula and Varadhan's lemma. Here we prove that these LD results are simple consequences of the LDP for
$(\mathcal{L}^n)_{n \geq 1}$. In this way we give a more direct and less technical proof of the LDPs in \cite{a06d,a06c}. 

Symmetrisation of sequences or vectors of random variables was already present in the literature prior to \cite{a06d,a06c}.
In \cite{dz92} the LD behavior of the empirical measures $\frac{1}{n} \sum_{i=1}^{n} \delta_{Y_i}$ of an exchangeable
sequence ${\bf Y} = (Y_1,\dots,Y_n,\dots)$ of random variables is investigated. The analysis relies on de Finetti's 
theorem which identifies the distribution of ${\bf Y}$ with the mixture of distributions of sequences of 
independent and identically distributed random variables. In \cite{t02} an LDP for the sequence of 
empirical measure processes $t \in [0,1] \mapsto \frac{1}{n} \sum_{i=1}^{[nt]} \delta_{X_i^n}$ 
is obtained where the $(X_i^n)_{1 \leq i \leq n}$ are symmetrised vectors in the sense that

$$
(X^n_1,\dots,X^n_n) = \frac{1}{n!} \sum_{\sigma \in \mathfrak{S}_n} (Y^{n}_{\sigma(1)},\dots,Y^{n}_{\sigma(n)})
$$
for a given triangular array $((Y_i^n)_{1 \leq i \leq n})_{n \geq 1}$ of random variables such that 
$\frac{1}{n} \sum_{i=1}^{n} \delta_{Y^n_i}$ satisfies an LDP.

\bigskip

\bigskip

Now let us describe our results more precisely. First we consider the following particular case: We 
are given a triangular array $((x_{i}^{n})_{1 \leq i \leq n})_{n \geq 1}$ of fixed elements of $\Sigma$, possibly 
with repetition, and we assume that 

$$
\mu^n = \frac{1}{n} \sum_{i=1}^{n} \delta_{x_{i}^{n}} \cvloi \mu \in M^1(\Sigma)
$$
where $\cvloi$ stands for weak convergence in the set $M^1(\Sigma)$ of Borel probability measures on $\Sigma$. 
Let $(V^n)_{n \geq 1}$ be the sequence of random measures defined on $(E,\mathcal{E},\mathbb{P})$ by

$$
V^n = \frac{1}{n} \sum_{i=1}^{n} \delta_{(x_i^n,x_{\sigma_n(i)}^n)} \in M^1(\Sigma^2).
$$
In order to describe the LD behavior of $(V^n)_{n \geq 1}$ we need to introduce some notations. For 
every $\nu \in M^1(\Sigma^2)$ we denote by $\nu_1$ (resp. $\nu_2$) its first (resp. second) marginal. 
For any two probabilities $\rho,\nu$ on a measurable space $(A,\mathcal{A})$ we denote by

$$
H(\nu| \rho) = \left\{
\begin{array}{cl}
\int_A {\mbox d}\nu \log \frac{{\mbox d}\nu}{{\mbox d}\rho} & \mbox{if \ } \nu << \rho \\
+ \infty & \mbox{otherwise}
\end{array}
\right.
$$
the relative entropy of $\nu$ with respect to $\rho$. Our first result is the following

\bigskip

\begin{theorem} The sequence $(V^n)_{n \geq 1}$ satisfies an LDP on $M^1(\Sigma^2)$ endowed with the weak convergence 
topology with good rate function

$$
\mathcal{I}(\nu) = \left\{
\begin{array}{cl}
H(\nu| \mu \otimes \mu) & \mbox{if \ } \nu_1 = \nu_2 = \mu \\
+ \infty & \mbox{otherwise}.
\end{array}
\right.
$$
\end{theorem}

\bigskip
\noindent
an LDP for $(\mathcal{L}^n)_{n \geq 1}$ follows easily

\bigskip

\begin{theorem}
Suppose that $(\mathcal{L}^n_1 = \frac{1}{n} \sum_{i=1}^{n} \delta_{X_{i}^{n}})_{n \geq 1}$ satisfies 
an LDP on $M^1(\Sigma)$ endowed with the weak convergence topology with good rate function $\mathcal{S}$. Then 
$(\mathcal{L}^n)_{n \geq 1}$ follows an LDP on $M^1(\Sigma^2)$
endowed with the weak convergence topology with good rate function

$$
\mathcal{J}(\nu) = \left\{
\begin{array}{cl}
\mathcal{S}(\nu_1) + H(\nu| \nu_1 \otimes \nu_1) & \mbox{if \ } \nu_1 = \nu_2 \\
+ \infty & \mbox{otherwise}.
\end{array}
\right.
$$
\end{theorem}

\bigskip

\noindent
Indeed, we can consider $\mathcal{L}^n$ as the result of a two steps mechanism: 
First we pick the $(X_i^n)_{1 \leq 1 \leq n}$ and, once these values are given, we pick the $\sigma_n$ and 
compute the resulting $\mathcal{L}^n$. While this may look artificial (since $\sigma_n$ is independent
from $(X_i^n)_{1 \leq 1 \leq n}$) it allows to identify the distribution of $\mathcal{L}^n$ as a mixture of
Large Deviations Systems (hereafter abbreviated LDS), see \cite{daga}. Then Theorem 2 follows directly from 
Theorem 1 thanks to Theorem 2.3 in \cite{g96}. Furthermore, by keeping in mind the mixture of LDS 
representation, it is an easy task to extend Theorem 2 to obtain LD results for a broad class of symmetrised bridge 
processes including those considered in \cite{a06d,a06c}.
For the sake of brevity we only give a statement on the LD behavior of the empirical path measure for continuous path
symmetrised bridge processes. It should be clear to the reader how other LD results in \cite{a06d,a06c} can be extended by 
following our approach.

More precisely, let $\mathbb{R}^d$ be a Polish space and $\mathcal{C} = \mathcal{C}([0,\beta],\mathbb{R}^d)$ be the Polish 
space of $\mathbb{R}^d$-valued continuous
functions defined on the interval $[0,\beta] \  (\beta > 0$ fixed) endowed with the uniform convergence metric. Let 
$\xi = (\xi_s)_{s \in [0,\beta]}$ be a 
$\mathcal{C}$-valued random variable defined on $(E,\mathcal{E},\mathbb{P})$. We assume that 
for every $(x,y) \in (\mathbb{R}^d)^2$ the distribution of $\xi$ conditioned 
on $\{\xi_0 = x, \xi_{\beta} = y \}$ is well-defined and denote it by $\mathbb{P}_{x,y}^{\xi}$. 
We further assume that for every element $\phi$ of the set $\mathcal{C}_b(\mathcal{C})$ 
of real valued bounded continuous function  defined on $\mathcal{C}$, the map

$$(x,y) \in (\mathbb{R}^d)^2 \mapsto \mathbb{E}_{x,y}^{\xi} e^{\phi(\xi)}$$
is bounded and continuous. Let $((X_{i}^{n})_{1 \leq i \leq n})_{n \geq 1}$ be a a triangular 
array of $\mathbb{R}^d$-valued random variables defined on $(E,\mathcal{E},\mathbb{P})$ such 
that $\frac{1}{n} \sum_{i=1}^{n} \delta_{X_{i}^{n}}$ 
satisfies an LDP on $M^1(\mathbb{R}^d)$ endowed with the weak convergence topology, with good rate function
$\mathcal{S}$. For every $n \geq 1$ we consider the $\mathcal{C}^n$-valued vector $(\xi^1,\dots,\xi^n)$ which 
distribution is given by 

$$\mathbb{P}^{sym}_n = \frac{1}{n!} \sum_{\sigma \in \mathfrak{S}_n} 
\int_{(\mathbb{R}^d)^n} \mathbb{P}(X^n_1 \in dx_1,\dots, X^n_n \in dx_n) \bigotimes_{i=1}^{n} 
\mathbb{P}^{\xi}_{x_i,x_{\sigma(i)}}$$
and the corresponding empirical path measure

\begin{equation}
L_n = \frac{1}{n} \sum_{i=1}^{n} \delta_{\xi^i}.
\label{pro}
\end{equation}
For every $\mu \in M^1(\mathcal{C})$ let us denote by $\mu_0$ (resp. $\mu_{\beta}$) the marginal distribution at time
0 (resp. at time $\beta$). By $\mu_{0,\beta}$ we denote the joint distribution at initial 
and final times 0 and $\beta$. Finally, we introduce the rate function $L$ defined on $M^1(\mathcal{C})$ by

$$
L(\mu) = \sup_{\phi \in \mathcal{C}_b(\mathcal{C})} \left\{ \langle \phi , \mu \rangle - \int_{(\mathbb{R}^d)^2}
\mu_{0,\beta}(dx,dy) \log \mathbb{E}^{\xi}_{x,y} e^{\phi(\xi)} \right\}.
$$
The following result is a consequence of Theorem 2 in the same way Theorem 2 is a consequence of Theorem 1.

\bigskip

\begin{theorem}
The sequence $(L_n)_{n \geq 1}$ obeys an LDP on $M^1(\mathcal{C})$ endowed with the weak convergence topology
with good rate function

$$
\mathcal{T}(\mu) = \left\{
\begin{array}{cl}
\mathcal{S}(\mu_0) + H(\mu_{0,\beta}| \mu_0 \otimes \mu_{\beta}) + L(\mu) & \mbox{if \ } \mu_0 = \mu_{\beta} \\
+ \infty & \mbox{otherwise}.
\end{array}
\right.
$$

\end{theorem}

\bigskip

\noindent
{\bf Remarks}\\
1) Theorem 1.1 in \cite{a06c} is a particular case of our Theorem 3 where $\xi$ is the Brownian motion and the 
$((X_{i}^{n})_{1 \leq i \leq n})_{n \geq 1}$ are independent and have the same 
distribution $\mathfrak{m}$. There it is further assumed that $\mathfrak{m}$ has compact support. Actually our result also holds
for general Polish spaces instead of $\mathbb{R}^d$.\\
2) Applications of Theorem 3 to mean-field interacting Bosons systems with a smooth coupling functional are possible.
For examples of such applications see \cite{a06d}.\\
3) Finally we would like to underscore the fact that there is a  3 step mechanism involved in the definition 
of the symmetrised bridge processes considered here. The conciseness and simplicity of our proof follows from this remark. 
This fact is also well illustrated by the particular form of the rate function $\mathcal{T}$.

\bigskip

\bigskip

There are three Sections in the remainder of the paper. Each one is devoted to the proof of one of our results.


\section{Proof of Theorem 1}


\noindent
Throughout this section we are given a triangular array $((x_{i}^{n})_{1 \leq i \leq n})_{n \geq 1}$ of elements of
$\Sigma$, possibly with repetition, and we assume that 

$$
\mu^n = \frac{1}{n} \sum_{i=1}^{n} \delta_{x_{i}^{n}} \cvloi \mu.
$$
For every $n \geq 1$ we note 

$$
\mathcal{V}_n = \left\{ \nu \in M^1(\Sigma^2) : \exists \sigma \in \mathfrak{S}_n,  \ 
\nu = \frac{1}{n} \sum_{i=1}^{n} \delta_{(x_{i}^{n}, x_{\sigma(i)}^{n})} \right\}.
$$
We introduce two triangular arrays of $\Sigma-$valued random variables 
$((L_{i}^ {n})_{1 \leq i \leq n})_{n \geq 1}$ and $((R_{i}^ {n})_{1 \leq i \leq n})_{n \geq 1}$ defined on a 
probability space $(\Omega, \mathcal{A}, \P)$ such that for every $n \geq 1$ the $2n$ random variables 
$L^{n}_{1},\dots, L^{n}_{n}, R^{n}_{1},\dots, R^{n}_{n}$ are mutually independent and each of them 
is distributed according to $\mu^n$. The sequence of random empirical measures

$$
W^n = \frac{1}{n} \sum_{i=1}^{n} \delta_{(L_{i}^{n}, R_{i}^{n})} \in M^1(\Sigma^2)
$$
defined on $(\Omega, \mathcal{A}, \P)$ has the following LD behavior


\begin{lemma}
The sequence $(W^n)_{n \geq 1}$ satisfies an LDP on $M^1(\Sigma^2)$ endowed with the weak convergence 
topology with good rate function $H(\nu | \mu \otimes \mu)$.
\label{lem1}
\end{lemma}


\noindent
{\bf Proof}
Since $\mu^{n} \cvloi \mu$ we have $\mu^{n}\otimes \mu^{n} \cvloi \mu \otimes \mu$ 
(see \cite{b68}, Chapter 1, Theorem 3.2). The announced result then follows from Theorem 3 
in \cite{bj88}. \hfill $\Box$

\bigskip

\noindent
Our strategy in proving Theorem 1 consists in comparing $V^n$ to random measures coupled to $W^n$. Comparison 
is possible because the $\nu \in M^1(\Sigma^2)$ such that $\mathcal{I}(\nu) < + \infty$ can 
be approached in the weak convergence topology by elements of $\mathcal{V}_n$. Our proof of this property
requires to use several metrics on $M^1(\Sigma^2)$ compatible with the weak convergence topology. 
This is the reason why in Section 2.1 we give a short account on the weak convergence topology 
prior to the proof of our approximation result. In Section 2.2 
we construct our coupling and finally in Section 2.3 we give the proper proof of Theorem 1.


\subsection{An approximation result}


We are given a Polish space $(\Sigma,d)$. The distance $d$ is a priori not a bounded metric but 
it is topologically equivalent to the bounded metric

$$
\tilde{d}(x,y) = \frac{d(x,y)}{1 + d(x,y)}
$$
which makes $\Sigma$ a separable metric space but not necessarily a complete one. The product topology on $\Sigma^2$ 
is metrizable by e.g. 

\begin{equation}
d_{2,M}((x_1,x_2),(y_1,y_2)) = \max ( d(x_1,y_1),d(x_2,y_2))
\label{dm}
\end{equation}
or 

\begin{equation}
d_{2,+}((x_1,x_2),(y_1,y_2)) = d(x_1,y_1) + d(x_2,y_2).
\label{d+}
\end{equation}
They both make $\Sigma^2$ a Polish space. We can also metrize the product topology on $\Sigma^2$ with the analogues 
$\tilde{d}_{2,M}$ and $\tilde{d}_{2,+}$ of (\ref{dm}) and (\ref{d+}) built on the ground of $\tilde{d}$.

The weak convergence topology on $M^1(\Sigma)$ is metrizable by the so-called dual-bounded-Lipschitz metric

\begin{equation}
\beta_{BL,\delta}(\rho,\nu) = 
\sup_{f \in  C_{b}(\Sigma) \atop \|f\|_{\infty} + \|f\|_{L,\delta} \leq 1} \left\{ \left| \int_{\Sigma} f d\rho - \int_{\Sigma} f d\nu \right| \right\}  \label{beta}
\end{equation}

\noindent
where $C_{b}(\Sigma)$ stands for the set of bounded continuous functions defined on $\Sigma$,

$$
  \|f\|_{\infty} = \sup_{x \in \Sigma} |f(x)| \ \ \mbox{and} \ \ 
  \|f\|_{L,\delta} = \sup_{x,y \in \Sigma \atop x \neq y} \left| \frac{f(x)-f(y)}{\delta(x,y)} \right|,
$$
where $\delta$ is either $d$ or $\tilde{d}$ (see \cite{d02}, Chapter 11, Theorem 11.3.3).
According to Kantorovitch-Rubinstein Theorem (see \cite{d02}, Chapter 11, Theorem 11.8.2)
the following metric on $M^1(\Sigma)$

$$
\beta_{W,\tilde{d}}(\rho,\nu) = \inf_{Q \in M^1(\Sigma^2) \atop 
Q_1 = \rho \ Q_2 = \nu} \left\{ \int_{\Sigma^2} \tilde{d}(x,y) Q(dx,dy)   \right\},
$$
the so-called Wasserstein metric associated to $\tilde{d}$ is compatible with the weak convergence topology as well.
However, note that the "analogue" of $\beta_{W,\tilde{d}}$ built on the ground of $d$ is in general not a 
metric for the weak convergence 
topology (for an illustration of this fact see \cite{d02} p.420-421). Finally we shall denote by 

\begin{equation}
\beta_{BL,\tilde{d}_{2,M}}(\rho,\nu) = 
\sup_{f \in  C_{b}(\Sigma^2) \atop \|f\|_{2,\infty} + \|f\|_{L,\tilde{d}_{2,M}} \leq 1} 
\left\{ \left| \int_{\Sigma^2} f d\rho - \int_{\Sigma^2} f d\nu \right| \right\}  \label{betap}
\end{equation}

\noindent
with 

$$
  \|f\|_{2,\infty} = \sup_{(x,y) \in \Sigma^2} |f(x,y)| \ \ \mbox{and} \ \ 
  \|f\|_{L,\tilde{d}_{2,M}} = \sup_{(x_1,y_1),(x_2,y_2) \in \Sigma^2, \atop (x_1,y_1) \neq (x_2,y_2)} \left| 
  \frac{f(x_1,y_1)-f(x_2,y_2)}{\tilde{d}_{2,M}((x_1,y_1),(x_2,y_2))} \right|,
$$
and  

$$
\beta_{W,\tilde{d}_{2,+}}(\rho,\nu) = \inf_{Q \in M^1(\Sigma^2 \times \Sigma^2) \atop Q_1 = \rho \ Q_2 = \nu} 
\left\{ \int_{\Sigma^2 \times \Sigma^2} \tilde{d}_{2,+}(x,y) Q(dx,dy) \right\},
$$
two metrics on $M^1(\Sigma^2)$ compatible with the weak convergence topology. The following lemma 
is a key result in the proof of Theorem 1. 


\begin{lemma}
Let $\nu \in M^1(\Sigma^2)$ be such that $\nu_1 = \nu_2 = \mu$. For every $n \geq 1$ there exists a
$\nu_n \in \mathcal{V}_n$ such that $\nu_n \cvloi \nu$.
\label{lem2}
\end{lemma}


\noindent
{\bf Proof}
Let $\nu \in M^1(\Sigma^2)$ be such that $\nu_1 = \nu_2 = \mu$. According to Varadarajan's Lemma (see \cite{d02}, Chapter 11,
Theorem 11.4.1) there exists a family $((u_{i}^{n},v_{i}^{n} )_{1 \leq i \leq n})_{n \geq 1}$ of 
elements of $\Sigma^2$ such that

$$
\gamma^n = \frac{1}{n} \sum_{i=1}^{n} \delta_{(u_{i}^{n}, v_{i}^{n})} \cvloi \nu.
$$
For every $n \geq 1$ we take $\sigma_n, \tau_n \in \mathfrak{S}_n$ such that

$$
\sum_{i=1}^{n} \tilde{d}(u_{i}^{n}, x_{\sigma_n(i)}^{n})
= \min_{\sigma \in \mathcal{V}_n} \left\{  \sum_{i=1}^{n} 
\tilde{d}(u_{i}^{n}, x_{\sigma(i)}^{n}) \right\}
$$

\noindent
and

$$
\sum_{i=1}^{n} \tilde{d}(v_{i}^{n}, x_{\tau_n(i)}^{n})
= \min_{\tau \in \mathcal{V}_n} \left\{  \sum_{i=1}^{n} 
\tilde{d}(v_{i}^{n}, x_{\tau(i)}^{n}) \right\}.
$$
We shall prove that 

$$
\nu^n = \frac{1}{n} \sum_{i=1}^{n} \delta_{(x_{\sigma_n(i)}^{n},x_{\tau_n(i)}^{n})} \cvloi \nu.
$$
To this end it is sufficient to prove that $\beta_{BL,\tilde{d}_{2,M}}(\nu^n,\nu) \rightarrow 0$. 
Let $\varepsilon > 0$ be fixed. There exists an $N_0$ such that for every $n \geq N_0$

\begin{equation}	
\beta_{BL, \tilde{d}_{2,M}}(\nu, \gamma^n) <  \varepsilon/3.
\label{pe}
\end{equation}

\noindent
Since $\gamma^n_1 \cvloi \mu$ there exists an $N_1$ such that for 
every $n \geq N_1$ we have $\beta_{W,\tilde{d}}(\gamma^n_1, \mu^n) < \varepsilon / 3$.
In a first step we will show that due to this for every $n \geq N_1$ we have

\begin{equation}
\frac{1}{n} \sum_{i=1}^{n} \tilde{d}(u_{i}^{n}, x_{\sigma_n(i)}^{n}) < \varepsilon / 3.
\label{se}
\end{equation}
In a second step we will show that (\ref{se}) leads to 

\begin{equation}
\beta_{BL, \tilde{d}_{2,M}}(\hat{\gamma}^n, \gamma^n) < \varepsilon / 3
\label{te}
\end{equation}
for every $n \geq N_1$ where

$$\hat{\gamma}^n = \frac{1}{n} \sum_{i=1}^{n} 
\delta_{(x_{\sigma_n(i)}^{n}, v_{i}^{n})}.$$ 
Analogously  to (\ref{te}) one can prove that there exists an $N_2$ such that for every $n \geq N_2$

\begin{equation}
\beta_{BL,\tilde{d}_{2,M}}( \hat{\gamma}^n, \nu^n) 
< \varepsilon / 3.
\label{qe}
\end{equation}
By combining (\ref{pe}, \ref{te}, \ref{qe}) we obtain the announced result.\\

\noindent
\underline{\it Step 1.} Since $\gamma_1^n$ and $\mu^n$ have finite support every Borel probability 
measure $Q$ on $\Sigma^2$ such that $Q_1 = \gamma_{1}^{n}$ and $Q_2 = \mu^n$ is of the form

\begin{equation}
Q(\alpha) = \frac{1}{n} \sum_{i,j=1}^{n} \alpha_{i,j} \delta_{(u_{i}^{n},x_{j}^{n})}
\label{bisto}
\end{equation}
where $\alpha = (\alpha_{i,j})_{1 \leq i,j \leq n}$ is an $n \times n$ bi-stochastic matrix.
Conversely every $n \times n$ bi-stochastic matrix $\alpha$ defines through (\ref{bisto}) 
a Borel probability measure $Q(\alpha)$ on $\Sigma^2$ such that $Q(\alpha)_1 = \gamma_{1}^{n}$ 
and $Q(\alpha)_2 = \mu^n$. Actually every bi-stochastic matrix can be written
as a convex combination of permutation matrices. These are $n \times n$ matrices 
with a single 1 in every line and every column, all other entries being 0 (for a proof of this 
fact see e.g. \cite{l68}, Chapter 11, Example 11.2). There is an obvious one-to-one correspondence 
between elements of $\mathfrak{S}_n$ and $n \times n$ permutation matrices. For 
every $\delta \in \mathfrak{S}_n$ we shall denote by $K_{\delta}$ the permutation matrix naturally associated to it. Therefore
for any Borel probability measure $Q$ on $\Sigma^2$ such that $Q_1 = \gamma_{1}^{n}$ and $Q_2 = \mu^n$
i.e. any choice of the components $(\lambda_{\delta})_{\delta \in \mathfrak{S}_n}$ of the convex 
combination $\alpha = \sum_{\delta \in \mathfrak{S}_n} \lambda_{\delta} K_{\delta}$ such that 
$Q = Q(\alpha)$ we have

\begin{eqnarray}
\int_{\Sigma^2} \tilde{d}(x,y) Q(dx,dy) & = & \sum_{\delta \in \mathfrak{S}_n} 
                                              \lambda_{\delta} (\sum_{i=1}^{n} \tilde{d}(u_{i}^{n}, x_{\delta(i)}^{n})) \\
                                        & \geq & \sum_{i=1}^{n} \tilde{d}(u_{i}^{n}, x_{\sigma_n(i)}^{n}).
\end{eqnarray}
Hence for every $n \geq N_1$

\begin{eqnarray*}
\frac{1}{n} \sum_{i=1}^{n} \tilde{d}(u_{i}^{n}, x_{\sigma_n(i)}^{n}) & = &
\beta_{W,\tilde{d}}( \gamma_1^n, \mu^n)
 \\
& < & \varepsilon / 3.
\end{eqnarray*}

\noindent
which proves (\ref{se}).\\

\noindent
\underline{\it Step 2.} For every $f \in \Cb$ such that 
$||f||_{2,\infty} + ||f||_{L,\tilde{d}_{2,M}} \leq 1$ and every $n \geq N_1$ we have

\begin{eqnarray*} 
\left| \int_{\Sigma^2} f d \hat{\gamma}^n - \int_{\Sigma^2} f d \gamma^n \right| & \leq & \frac{1}{n} \sum_{i=1}^{n}
\left|f(x_{\sigma_n(i)}^{n}, v_{i}^{n}) - f(u_{i}^{n}, v_{i}^{n}) \right| \\
& \leq & \frac{1}{n} \sum_{i=1}^{n} \tilde{d}_{2,M} ((x_{\sigma_n(i)}^{n}, v_{i}^{n}), (u_{i}^{n}, v_{i}^{n})) \\ 
& \leq & \frac{1}{n} \sum_{i=1}^{n} \tilde{d}(u_{i}^{n}, x_{\sigma_n(i)}^{n}) \\
& < & \varepsilon / 3
\end{eqnarray*}

\noindent
hence $\beta_{BL, \tilde{d}_{2,M}}(\hat{\gamma}^n, \gamma^n) < \varepsilon / 3$ for every $n \geq N_1$.  \hfill $\Box$ 


\subsection{Coupled empirical measures}


To every $n \geq 1$ and every realization of $W^n$ we associate two elements of $M^1(\Sigma^2)$ by 

\begin{equation} 
\beta_{W,\tilde{d}_{2,+}}(W^n,\widetilde{W}^n) = 
\min_{\nu \in \mathcal{V}_n} \left\{ \beta_{W,\tilde{d}_{2,+}}(W^n,\nu) \right\},
\label{argmin}
\end{equation}

\noindent
and  

\begin{equation} 
\beta_{W,\tilde{d}_{2,+}}(W^n,\widehat{W}^n) = 
\max_{\nu \in \mathcal{V}_n} \left\{ \beta_{W,\tilde{d}_{2,+}}(W^n,\nu) \right\},
\label{argmax}
\end{equation}

\noindent
In case there are several elements of $\mathcal{V}_n$ achieving the 
min (resp. the max) $\widetilde{W}^n$ (resp. $\widehat{W}^n$) is 
picked uniformly at random among these measures.

\begin{lemma} For every $n \geq 1$ the random measures $\widetilde{W}^n, \widehat{W}^n$ and $V^n$ are 
identically distributed over $M^1(\Sigma^2)$. 
\label{lem3}
\end{lemma}

\noindent
{\bf Proof}
We shall only prove that $\widetilde{W}^n$ and $V^n$ are identically distributed 
since the proof with $\widehat{W}^n$ and $V^n$
works the same way. Let $n \geq 1$ be fixed. For the sake of clarity let us first assume that there is no repetition among the
$x_1^n,\dots,x_n^n$. In this case every $\nu \in \mathcal{V}_n$ corresponds to a single $\tau \in \mathfrak{S}_n$ by

\begin{equation}
\nu = \frac{1}{n} \sum_{i=1}^{n} \delta_{(x_{i}^{n}, x_{\tau(i)}^{n})}.
\label{code}
\end{equation}
Next let us consider a fixed realization $(l_i^n,r_i^n)_{1 \leq i \leq n}$ of $(L_i^n,R_i^n)_{1 \leq i \leq n}$
and let us denote by $w^n$ the corresponding value for $W^n$. Due to the property of the minimizer in the 
Wasserstein distance between two atomic measures we already employed in the proof of Lemma 1, for 
every $\nu \in \mathcal{V}_n$ (i.e. every $\tau \in \mathfrak{S}_n$ according to (\ref{code})) there 
exists a $\sigma \in \mathfrak{S}_n$ such that
\begin{eqnarray*}
\beta_{W,\tilde{d}_{2,+}}(w^n,\nu) & = & \frac{1}{n} \sum_{i=1}^{n} \tilde{d}_{2,+} ((l_i^n,r_i^n), 
(x_{\sigma(i)}^{n}, x_{\sigma \circ \tau(i)}^{n})) \\
 & = & \frac{1}{n} \sum_{i=1}^{n} \tilde{d}(l_i^n,x_{\sigma(i)}^{n}) 
+ \frac{1}{n} \sum_{i=1}^{n} \tilde{d}(r_i^n,x_{\sigma \circ \tau(i)}^{n}).
\end{eqnarray*}

\noindent
Thus, since we are looking for the maximum over $\sigma$ and $\tau$, for a fixed 
realization $(l_i^n,r_i^n)_{1 \leq i \leq n}$ of $(L_i^n,R_i^n)_{1 \leq i \leq n}$ the 
corresponding value of $\widetilde{w}^n$ is obtained by finding $\eta_l$ and $\eta_r$ such that 

$$
\sum_{i=1}^{n} \tilde{d} (l_i^n,x_{\eta_l(i)}^n) = \min_{\sigma \in \mathfrak{S}_n} \left\{ \sum_{i=1}^{n} 
\tilde{d} (l_i^n,x_{\sigma(i)}^n) \right\}
$$
and
$$
\sum_{i=1}^{n} \tilde{d} (r_i^n,x_{\eta_r(i)}^n) = \min_{\sigma \in \mathfrak{S}_n} \left\{ \sum_{i=1}^{n} 
\tilde{d} (r_i^n,x_{\sigma(i)}^n) \right\}
$$
and taking

$$
\tilde{w}^n = \frac{1}{n} \sum_{i=1}^{n} \delta_{(x_{\eta_l(i)}^n,x_{\eta_r(i)}^n)}.
$$
In case several $\eta_l$ and/or $\eta_r$ realise the minima in the displays above, those defining $\widetilde{w}^n$ 
are picked among them uniformly at random. Now, remark that for every $\gamma_{l},\gamma_{r} \in \mathfrak{S}_n$,  
observing $(l_{\gamma_l(i)}^{n}, r_{\gamma_r(i)}^{n})_{1 \leq i \leq n}$ has the same probability as observing  
$(l_i^n,r_i^n)_{1 \leq i \leq n}$ and results in $\gamma_l \circ \eta_l$ and $\gamma_r \circ \eta_r$ in defining
$\widetilde{w}^n$ instead of $\eta_l$ and $\eta_r$. Thus if we consider $\eta_l$ and $\eta_r$ as random variables 
defining $\widetilde{W}^n$, we see that their distribution {\it conditioned on }$W^n$ is uniform over $\mathfrak{S}_n$.
Hence $\widetilde{W}^n$ and $V^n$ are both uniformly distributed over $M^1(\Sigma^2)$, thus identically distributed. 
This proof extends easily to the case when there are repetitions among the $x_1^n,\dots,x_n^n$. \hfill $\Box$


\subsection{Proof of the LD bounds}


We start the proof of the LD bounds by proving the following

\begin{lemma} We have:
\begin{enumerate}
\item $\mathcal{I}$ is a good rate function.
\item The sequence $(V^n)_{n \geq 1}$ is exponentially tight. 
\end{enumerate}
\label{lem4}
\end{lemma}

\noindent
{\bf Proof}\\
(1) Let $\alpha \geq 0$. We have

\begin{eqnarray*}
N_{\alpha} & = & \{ \nu \in M^1(\Sigma^2) : \mathcal{I}(\nu) \leq \alpha \} \\
& = & \{ \nu \in M^1(\Sigma^2) : H(\nu | \mu \otimes \mu) \leq \alpha \} \cap 
\{ \nu \in M^1(\Sigma^2) : \nu_1 = \nu_2 = \mu \}.
\end{eqnarray*}

\noindent
Thus, for every $\alpha \geq 0$, $N_{\alpha}$ is the intersection of a compact and a closed subset of 
$M^1(\Sigma^2)$, therefore it is compact.\\

\noindent
(2) For every measurable $A \subset M^1(\Sigma^2)$ we have 

\begin{eqnarray*}
\limsup_{n \rightarrow \infty} \frac{1}{n} \log \mathbb{P}( V^n \in A^c) & = & \limsup_{n \rightarrow \infty} 
\frac{1}{n} \log \P( W^n \in A^c | \frac{1}{n} \sum_{i=1}^{n} 
\delta_{R_i^n} = \frac{1}{n} \sum_{i=1}^{n} \delta_{L_i^n} = \mu^n ) \\
& \leq & \limsup_{n \rightarrow \infty} \frac{1}{n} \log \P( W^n \in A^c) \\ 
& & \ \ \ \ \ \ \ \ - 2 \liminf_{n \rightarrow \infty} \frac{1}{n} 
\log \P(\frac{1}{n} \sum_{i=1}^{n} \delta_{L_i^n} = \mu^n )
\end{eqnarray*}

\noindent
Since $(W^n)_{n \geq 1}$ satisfies an LDP on $M^1(\Sigma^2)$ with a good rate function it is 
exponentially tight (see \cite{dz98}, Remark a) p.8). Thus for every $\alpha \geq 0$ we can chose a 
compact set $A_{\alpha} \subset M^1(\Sigma^2)$ that makes the first term in the last display 
smaller than $ - \alpha - 2$. On the other hand it is clear that
for every $n \geq 1$ we have

$$
\P( \frac{1}{n} \sum_{i=1}^{n} \delta_{L_i^n} = \mu^n ) \geq n!  \frac{1}{n^n}
$$
equality corresponding to the case when there are no ties among the $x^n_1,\dots, x^n_n$. Thus

$$
- 2 \liminf_{n \rightarrow \infty} \frac{1}{n} \log \P(\frac{1}{n} \sum_{i=1}^{n} \delta_{L_i^n} = \mu^n ) \leq 2 
$$
which completes the proof. \hfill $\Box$

\subsubsection{Proof of the lower bound}

It is sufficient in order to prove the lower bound of the LDP to prove that 

$$
- \mathcal{I}(\nu) \leq \liminf_{n \rightarrow \infty} \frac{1}{n} \log \mathbb{P}( V^n \in B(\nu,\varepsilon))
$$
holds for every $\nu \in M^1(\Sigma^2)$ and every $\varepsilon > 0$, where $B(\nu,\varepsilon)$ stands for 
the open ball centered at $\nu \in M^1(\Sigma^2)$ of radius $\varepsilon > 0$ for the $\beta_{W,\tilde{d}_{2,+}}$ 
metric. So let $\varepsilon > 0$ and 
$\nu \in M^1(\Sigma^2)$ be such that $\mathcal{I}(\nu) < + \infty$. In particular $\nu_1 = \nu_2 = \mu$. 
According to Lemma \ref{lem2} there exists a sequence $(\nu^n)_{n \geq 1}$ of elements of $M^1(\Sigma^2)$ 
such that $\nu^n \in \mathcal{V}_n$ and $\nu^n \cvloi \nu$. According to Lemma \ref{lem3} we have

\begin{eqnarray*}
\mathbb{P}(V^n \in B(\nu,\varepsilon)) & = & \P( \widetilde{W}^n \in B(\nu,\varepsilon)) \\
& \geq & \P( \beta_{W,\tilde{d}_{2,+}}(\widetilde{W^n},W^n) < \frac{\varepsilon}{3},  
\beta_{W,\tilde{d}_{2,+}}(W^n,\nu^n) < \frac{\varepsilon}{3}, 
\beta_{W,\tilde{d}_{2,+}}(\nu^n,\nu) < \frac{\varepsilon}{6}) \\
& \geq & \P( \beta_{W,\tilde{d}_{2,+}}(\nu^n,W^n) < \frac{\varepsilon}{3}, 
\beta_{W,\tilde{d}_{2,+}}(\nu^n,\nu) < \frac{\varepsilon}{6}) \\
\end{eqnarray*}

\noindent
since it follows from the definition of $\widetilde{W^n}$ that for every $\nu^n \in \mathcal{V}_n$ we get

$$
\beta_{W,\tilde{d}_{2,+}}(\widetilde{W^n},W^n) \leq \beta_{W,\tilde{d}_{2,+}}(\nu^n,W^n).
$$
On the other hand since $\nu^n \cvloi \nu$ we see that for $n$ large enough 
$\left\{ \beta_{W,\tilde{d}_{2,+}}(\nu^n,\nu) < \frac{\varepsilon}{6} \right\} = \Omega $. Thus, for those $n$'s

\begin{eqnarray*}
\P( \beta_{W,\tilde{d}_{2,+}}(\nu^n,W^n) < \frac{\varepsilon}{3}, 
\beta_{W,\tilde{d}_{2,+}}(\nu^n,\nu) < \frac{\varepsilon}{6}) & \geq & 
\P( \beta_{W,\tilde{d}_{2,+}}(\nu,W^n) < \frac{\varepsilon}{6}, 
\beta_{W,\tilde{d}_{2,+}}(\nu^n,\nu) < \frac{\varepsilon}{6}) \\
& \geq & \P( \beta_{W,\tilde{d}_{2,+}}(W^n,\nu) < \frac{\varepsilon}{6})
\end{eqnarray*}

\noindent
Finally, it follows from Lemma \ref{lem1} that
\begin{eqnarray*}
\liminf_{n \rightarrow \infty} \frac{1}{n} \log \mathbb{P}( V^n \in B(\nu,\varepsilon)) & \geq & 
\liminf_{n \rightarrow \infty} \frac{1}{n} 
\log \P( \beta_{W,\tilde{d}_{2,+}}(\nu,W^n) 
< \frac{ \varepsilon}{6}) \\
& \geq & - H(\nu | \mu \otimes \mu) = - \mathcal{I}(\nu).
\end{eqnarray*}

\subsubsection{Proof of the upper bound}

In order to prove the upper bound of the LDP, it is sufficient to prove that it holds for compact subsets of
$M^1(\Sigma^2)$. Indeed, since $(V^n)_{n \geq 1}$ is an exponentially tight sequence (see Lemma \ref{lem4}) 
the full upper bound will follow from Lemma 1.2.18 in \cite{dz98}. Let $A$ be a compact subset 
of $M^1(\Sigma^2)$ and let us denote by

$$
A_{\mu} = \left\{ \nu \in A  : \nu_1 = \nu_2 = \mu \right\}
$$
which is a compact subset of $M^1(\Sigma^2)$ as well. Since the weak convergence 
topology on $M^1(\Sigma^2)$ is compatible with the $\beta_{W,\tilde{d}_{2,+}}$ metric, it makes $M^1(\Sigma^2)$ 
a regular topological space: For every $\nu \in A$ such that $\nu \in A_{\mu}^c$ there exists $\varepsilon_{\nu} > 0$ such 
that $B(\nu, 2 \varepsilon_{\nu})  \cap A_{\mu} = \emptyset$. In particular
$\bar{B}(\nu, \varepsilon_{\nu})  \cap A_{\mu} = \emptyset$ where $\bar{B}(\nu, \varepsilon)$ denotes the
closed ball centered on $\nu \in M^1(\Sigma^2)$ of radius $\varepsilon > 0$ for the $\beta_{W,\tilde{d}_{2,+}}$ metric. 
On the other hand, since $\nu \mapsto H(\nu | \mu \otimes \mu)$ is lower semi-continuous, for 
every $\nu \in A_{\mu}$ and every $\delta > 0$ there exists a $\varphi(\nu, \delta) > 0$ such that

$$
\inf_{\rho \in \bar{B}(\nu, \varphi(\nu, \delta))} H(\rho | \mu \otimes \mu) \geq (H(\nu | \mu \otimes \mu) - \delta) 
\wedge \frac{1}{\delta}.
$$
For every $\delta > 0$ we consider the coverage 

$$
A \subset \left( \cup_{\nu \in A \cap A_{\mu}^c} B (\nu, \varepsilon_{\nu}) \right) \cup \left( \cup_{\nu \in  A_{\mu}} 
B (\nu, \frac{\varphi(\nu, \delta)}{8}) \right) 
$$
from which we extract a finite coverage

$$
A \subset \left( \cup_{\nu \in I_1} B (\nu, \varepsilon_{\nu}) \right) \cup \left( \cup_{\nu \in I_2} 
B (\nu, \frac{\varphi(\nu, \delta)}{8}) \right) 
$$
where $I_1 \subset A \cap A_{\mu}^c$ and $I_2 \subset A_{\mu}$ are finite sets. Then, according to Lemma 1.2.15
in \cite{dz98} 

\begin{eqnarray*} 
\limsup_{n \rightarrow \infty} \frac{1}{n} \log \mathbb{P}( V^n \in A) & \leq & \max 
\left\{ \limsup_{n \rightarrow \infty} \frac{1}{n} \log \mathbb{P}( V^n \in  
\cup_{\nu \in I_1} \bar{B} (\nu, \varepsilon_{\nu}) \cap A), \right. \\
& & \ \ \ \ \left. \limsup_{n \rightarrow \infty} \frac{1}{n} \log \mathbb{P}( V^n \in  \cup_{\nu \in I_2} 
B (\nu, \frac{\varphi(\nu, \delta)}{8})) \right\}.
\end{eqnarray*}

\noindent
For every $\nu \in I_1$ there can not be an infinite number of elements of $\cup_{n \geq 1} \mathcal{V}_n$ in 
$\bar{B} (\nu, \varepsilon_{\nu}) \cap A$ for otherwise we would get $\bar{B} (\nu, \varepsilon_{\nu}) \cap A_{\mu} \neq \emptyset$.
The first term in the max is then equal to $- \infty$. We are left with the second term and according to Lemmas 1,2 and 3
we have

\begin{eqnarray*}
\limsup_{n \rightarrow \infty} \frac{1}{n} \log \mathbb{P}( V^n \in A) & \leq & 
\limsup_{n \rightarrow \infty} \frac{1}{n} \log 
\mathbb{P}( V^n \in  \cup_{\nu \in I_2} B (\nu, \frac{\varphi(\nu, \delta)}{8})) \\
& \leq & 
\max_{\nu \in I_2} \left\{ \limsup_{n \rightarrow \infty} \frac{1}{n} \log \P( \widehat{W}^n \in B (\nu, \frac{\varphi(\nu, \delta)}{8})) 
\right\} \\
& \leq & \max_{\nu \in I_2} \left\{ \limsup_{n \rightarrow \infty} \frac{1}{n} \log \P( \beta_{W,\tilde{d}_{2,+}}(\nu^n,\widehat{W}^n) 
< \frac{\varphi(\nu, \delta)}{4}) \right\} \\
& \leq & \max_{\nu \in I_2} \left\{ \limsup_{n \rightarrow \infty} \frac{1}{n} \log \P( \beta_{W,\tilde{d}_{2,+}}(\nu, W^n) 
< \frac{\varphi(\nu, \delta)}{2}) \right\} \\
& \leq & \max_{\nu \in I_2} \left\{ - \inf_{\rho \in \bar{B}(\nu, \varphi(\nu, \delta))} H(\rho | \mu \otimes \mu) \right\}\\
& \leq & \max_{\nu \in I_2} \left\{ - (H(\nu | \mu \otimes \mu) - \delta) \wedge \frac{1}{\delta} \right\} \\
& \leq & \max_{\nu \in I_2} \left\{ - (\mathcal{I}(\nu) - \delta) \wedge \frac{1}{\delta} \right\} \\
& \leq & - \inf_{\nu \in A} \left\{ (\mathcal{I}(\nu) - \delta) \wedge \frac{1}{\delta} \right\}.
\end{eqnarray*}

\noindent
By letting $\delta \rightarrow 0$ we obtain the announced upper bound, see Remark 1.2.10 in \cite{dz98}.


\section{Proof of Theorem 2}


\noindent
As mentioned in the Introduction, in order to prove Theorem 2 it is sufficient to prove that the distribution
of 

$$\mathcal{L}^n = \frac{1}{n} \sum_{i=1}^{n} \delta_{(X^n_{i},X^n_{\sigma_n(i)})}$$
on $M^1(\Sigma^2)$ is a mixture of LDS in the sense of \cite{g96}. For the sake of 
clarity we recover the notations of \cite{g96} when identifying the components of the LDS:

\bigskip

\noindent
$\bullet$ $\mathcal{Z} = M^1(\Sigma^2)$ is a Polish space when endowed with the weak convergence topology.\\
$\bullet$ $\mathcal{X} = M^1(\Sigma) (= \mathcal{X}_{\infty})$ is a Polish space when endowed with 
      the weak convergence topology as well.\\
$\bullet$ For every $n \geq 1$ we note
$$\mathcal{X}_n = \left\{ \nu \in M^1(\Sigma) : \exists (x_1,\dots,x_n) \in \Sigma^n \ \nu = \frac{1}{n}
\sum_{i=1}^{n} \delta_{x_i} \right\},$$
and according to Varadarajan's lemma for every $\nu \in \mathcal{X}$ there exists a sequence 
$(\nu^n)_{n \geq 1}$ such that for every $n \geq 1$ we have $\nu^n \in \mathcal{X}_n$ and $\nu^n \cvloi \nu$.\\
 $\bullet$ The map $\pi : \mathcal{Z} \rightarrow \mathcal{X}$ defined by $\pi(\nu) = \nu_1$ is continuous 
and surjective.\\
$\bullet$ For every $n \geq 1$ and every $\nu = \frac{1}{n}
\sum_{i=1}^{n} \delta_{x_i} \in \mathcal{X}_n$, let $P^n_{\nu}$ be the distribution of 
$V^n = \frac{1}{n} \sum_{i=1}^{n} \delta_{(x_i,x_{\sigma_n(i)})}$ under $\mathbb{P}$. The family 
$\Pi = \left\{ P^n_{\nu}, \nu \in  \mathcal{X}_n, n \geq 1 \right\}$
of finite measures on the Borel $\sigma-$field on $\mathcal{Z}$ is such that for every 
$n \geq 1$ and every $\nu \in \mathcal{X}_n$ we have $P^n_{\nu} (\pi^{-1}\{\nu\}^c) = 0.$\\
$\bullet$ Let $Q^n$ be the distribution of $\frac{1}{n} \sum_{i=1}^{n} \delta_{X^n_i}.$ For every $n \geq 1$ and every
measurable $A \subset M^1(\Sigma^2)$

$$\mathbb{P}(\mathcal{L}^n \in A) = \int_{\mathcal{X}_n} P^n_{\nu} (A) Q^n (d \nu).$$

\noindent
All the requirements of Definition 2.1 in \cite{g96} are satisfied by our model thanks to Theorem 1 .
It follows from Theorem 2.3 in \cite{g96} that the sequence $(\mathcal{L}^n)_{n \geq 1}$ obeys 
an LDP on $M^1(\Sigma^2)$ endowed with the weak convergence topology with good rate function

$$
\mathcal{J}(\nu) = \left\{
\begin{array}{cl}
\mathcal{S}(\nu_1) + H(\nu| \nu_1 \otimes \nu_1) & \mbox{if \ } \nu_1 = \nu_2 \\
+ \infty & \mbox{otherwise}.
\end{array}
\right.
$$


\section{Proof of Theorem 3}


\noindent
In order to prove Theorem 3 it is sufficient to prove that the distribution of 

$$
L_n = \frac{1}{n} \sum_{i=1}^{n} \delta_{\xi^i}
$$
under 

$$\mathbb{P}^{sym}_n = \frac{1}{n!} \sum_{\sigma \in \mathfrak{S}_n} 
\int_{(\mathbb{R}^d)^n} \mathbb{P}(X^n_1 \in dx_1,\dots, X^n_n \in dx_n) \bigotimes_{i=1}^{n} 
\mathbb{P}^{\xi}_{x_i,x_{\sigma(i)}}.$$ 
is a mixture of LDS. Once again, we use the notations of \cite{g96} to identify the components of the LDS:

\bigskip

\noindent
$\bullet$ $\mathcal{Z} = M^1(\mathcal{C})$ is a Polish space when endowed with the weak convergence topology.\\
$\bullet$ $\mathcal{X} = M^1(\mathbb{R}^d \times \mathbb{R}^d) (= \mathcal{X}_{\infty})$ is a Polish 
space when endowed with the weak convergence topology as well.\\
$\bullet$ For every $n \geq 1$ we note
$$\mathcal{X}_n = \left\{ \nu \in M^1(\mathbb{R}^d \times \mathbb{R}^d) : \exists ((s_1,a_1),\dots,(s_n,a_n)) 
\in (\mathbb{R}^d \times \mathbb{R}^d)^n \ \nu = \frac{1}{n} \sum_{i=1}^{n} \delta_{(s_i,a_i)} \right\},$$
and for every $\nu \in \mathcal{X}$ there exists a sequence 
$(\nu^n)_{n \geq 1}$ such that for every $n \geq 1$ we have $\nu^n \in \mathcal{X}_n$ and $\nu^n \cvloi \nu$.\\
$\bullet$ The map $\pi : \mathcal{Z} \rightarrow \mathcal{X}$ defined by $\pi(\nu) = \nu_{0,\beta}$ is continuous 
and surjective.\\
$\bullet$ For every $n \geq 1$ and every $\nu = \frac{1}{n}
\sum_{i=1}^{n} \delta_{(s_i,a_i)} \in \mathcal{X}_n$, let $P^n_{\nu}$ be the distribution of 
$L^n = \frac{1}{n} \sum_{i=1}^{n} \delta_{\xi^i}$ under $\bigotimes_{i=1}^{n} \mathbb{P}_{s_i,a_i}$.
The family $\Pi = \left\{ P^n_{\nu}, \nu \in  \mathcal{X}_n, n \geq 1 \right\}$
of finite measures on the Borel $\sigma-$field on $\mathcal{Z}$ is such that for every 
$n \geq 1$ and every $\nu \in \mathcal{X}_n$ we have $P^n_{\nu} (\pi^{-1}\{\nu\}^c) = 0.$\\
$\bullet$ Let $Q^n$ is the distribution of $\frac{1}{n} \sum_{i=1}^{n} \delta_{(X^n_i, X^n_{\sigma_n(i)})}.$ For every 
$n \geq 1$ and every measurable $A \subset M^1(\mathcal{C})$ 

$$\mathbb{P}^{sym}_n(\mathcal{L}^n \in A) = \int_{\mathcal{X}_n} P^n_{\nu} (A) Q^n (d \nu).$$
Thus, in order to prove Theorem 3, we are left to verify that for any $\mu \in M^1((\mathbb{R}^d)^2)$ and any
sequence $\mu^n = \frac{1}{n} \sum_{i=1}^{n} \delta_{(s_i,a_i)} \in \mathcal{X}_n$ such that $\mu^n \cvloi \mu$, the 
distribution of $L_n = \frac{1}{n} \sum_{i=1}^{n}\delta_{\xi^i}$ under $\bigotimes_{i=1}^{n} \mathbb{P}^{\xi}_{s_i,a_i}$
satisfies an LDP with good rate function $L$. For every $\phi \in \mathcal{C}_b(\mathcal{C})$ we have 

\begin{eqnarray*}
\Lambda_n(\phi) & = & \log \int_{\mathcal{C}^n} e^{n \langle \phi, L_n \rangle} \bigotimes_{i=1}^{n} \mathbb{P}^{\xi}_{s_i,a_i} \\ 
                & = & \sum_{i=1}^{n} \log \mathbb{E}^{\xi}_{s_i,a_i} e^{\phi(\xi)}
\end{eqnarray*}              
Thus 

\begin{eqnarray*}
\Lambda(\phi) & = & \lim_{n \rightarrow \infty} \frac{1}{n} \Lambda_n(\phi)\\
              & = & \int_{(\mathbb{R}^d)^2} \mu(dx,dy) \log \mathbb{E}^{\xi}_{x,y} e^{\phi(\xi)}.
\end{eqnarray*} 
The map $\phi \in \mathcal{C}_b(\mathcal{C}) \mapsto \Lambda(\phi)$ is G\^ateaux differentiable. On the other hand
the distribution of $L_n$ under $\bigotimes_{i=1}^{n} \mathbb{P}^{\xi}_{s_i,a_i}$ is exponentially tight as can be showed 
by adapting the proof of Lemma \ref{lem4}. We conclude thanks to Corollary 4.6.14 in \cite{dz98}.

\bigskip

\bigskip

\noindent
{\bf Acknowledgment} I thank Omer Adelman, Halim Doss and J\'er\^ome Renault for their interest 
in part of this work. I also thank the referee for his careful reading of the manuscript.

\bibliographystyle{plain}
\bibliography{bibfile}

\end{document}